\documentclass[final,3p,times]{elsarticle}

\usepackage{amssymb}
\usepackage{graphicx,setspace}
\usepackage{psfrag}
\usepackage{amsfonts}
\expandafter\let\csname equation*\endcsname\relax
\expandafter\let\csname endequation*\endcsname\relax
\usepackage{amsmath}
\usepackage{amssymb, amsthm}
\usepackage{mathrsfs}
\usepackage{epstopdf}
\usepackage{float}
\usepackage{lineno,hyperref}
\usepackage{color}
\usepackage{subfigure}

\usepackage{amsmath}
\usepackage{dsfont} 
\usepackage{amssymb} 
\usepackage{graphicx,color}
\usepackage{extarrows}

\usepackage{graphicx}

\usepackage{epstopdf}

\journal{arXiv}

\begin{document}
\newtheorem{definition}{Definition}[section]
\newtheorem{lemma}{Lemma}[section]
\newtheorem{remark}{Remark}[section]
\newtheorem{theorem}{Theorem}[section]
\newtheorem{proposition}{Proposition}
\newtheorem{assumption}{Assumption}
\newtheorem{example}{Example}
\newtheorem{corollary}{Corollary}[section]
\def\ep{\varepsilon}
\def\Rn{\mathbb{R}^{n}}
\def\Rm{\mathbb{R}^{m}}
\def\E{\mathbb{E}}
\def\hte{\hat\theta}
\renewcommand{\theequation}{\thesection.\arabic{equation}}
\begin{frontmatter}



\title{On the solitary wave configurations of nonlinear Schr\"odinger equation under the effect of L\'{e}vy noise}

\author{Hina Zulfiqar\fnref{addr1}}\ead{zhinazulfiqar@uo.edu.pk}
\author{Shenglan Yuan\corref{cor1}\fnref{addr2,addr3}}\ead{shenglanyuan@gbu.edu.cn}\cortext[cor1]{Corresponding author}
\author{Muhammad Shoaib Saleem\fnref{addr1}}\ead{shoaib83455@gmail.com}

\address[addr1]{\rm Department of Mathematics, University of Okara 56300, Pakistan}
\address[addr2]{\rm Department of Mathematics, School of Sciences, Great Bay University, Dongguan 523000, China}
\address[addr3]{\rm Great Bay Institute for Advanced Study, Songshan Lake International Innovation
Entrepreneurship Community A5, Dongguan 523000, China}
\begin{abstract}
This study aims to examine the effect of L\'{e}vy noise on the solutions of the nonlinear Schr\"odinger equation. An improved diversity of stochastic solutions is instinctively located discretely on certain conditions by applying the generalized Kudryashov method. Moreover, the dynamical behaviors of these exact results of the nonlinear Schr\"odinger equation are interpreted in the context of the effect of L\'{e}vy noise. Even mathematical evaluations have been conducted and presented.
\end{abstract}


\end{frontmatter}


\section{Introduction}
The nonlinear Schr\"odinger equation (NLS) exhibits an essential pattern for depicting wave structures that are found in several areas of science. The complexity of obtaining new solutions for nonlinear differential equations has been stimulating challenges for many centuries \cite{1,2}. It is intimate that nonlinear approaches are especially imperative in a multitude of scientific fields, notably nonlinear fiber optics, chemical physics, quantum mechanics, geochemistry, plasma waves, capillary gravity waves, biology, fluid dynamics and plasma physics \cite{3,4,5,6,7}. The NLS has a great abundance of applications, which is crucial for the description of nonlinear dissipation in various fields of basic science \cite{8,9,10}.

Consequently, numerous researchers focus on treating the equations of nonlinear processes to demonstrate their significance \cite{11,12,13,14,15}.  The NLS can be used to describe the propagation of optical beams in nonlinear media, such as optical fibers or crystals. The equation can also be applied to the study of ion-acoustic waves in plasmas. The Gross-Pitaevskii equation, which is a type of NLS, is used to describe the behavior of Bose-Einstein condensates. Furthermore, the NLS can be employed to model the propagation of water waves in certain conditions.

The specific consequences of random interruptions, thermal deviations, and spontaneous diffractions of nonlinear fluctuations were explored  by adapting stochastic phenomena to the Schr\"odinger equations \cite{16,17,18,19,20}. The existence and uniqueness of martingale solutions to a stochastic NLS driven by a L\'evy motion were established by Petroni and Pusterla \cite{21}. The fractional NLS was discussed by Raju, Panigrahi and Porsezian \cite{22}. The deliberated results for symmetric and antisymmetric solitons also originated from their work \cite{23}. The impact of the L\'evy indicator on distinct solitons was analyzed by Bo, Wang, Fang, \emph{et al}. \cite{24}. The reliability and solidity intervals of solitons were characterized by Cao and Dai \cite{25}. By using modulation or amplitude equations for the nonlinear stochastic partial differential equation influenced by cylindrical stable L\'evy processes, the approximation was established by Yuan and Bl\"{o}mker \cite{YB}.
The anti-interference potential of robust solitons against minor distraction exhibited adequate strength \cite{26,27,28}. Nevertheless, the classification of these equations has been subject to rigorous scrutiny in obtaining estimated or precise solutions \cite{29,30,31}. The fluctuations of some nonlinear partial differential equations (PDEs) with discrete types of noises were examined by many scholars \cite{32,33,34,35}.

The stochastic NLS with L\'{e}vy noise can be expressed as:
\begin{eqnarray}\label{1}
i \psi_{t}-\psi _{xx}+2|\psi|^2\psi-2\rho^2\psi+\sigma\psi L_{t} = 0, \hspace{0.5cm}\text{for}\hspace{0.2cm}t\geq0\hspace{0.3cm}\text{and}\hspace{0.2cm}x\in\mathbb{R}.
\end{eqnarray}
 Here, $\psi:=\psi(x,t)\in\mathbb{C}$ is a complex-valued function, $\rho$ and $\sigma$ are parameters, and $\frac{dL}{dt}$ = $L _{t}$ represents the first time derivative of L\'{e}vy motion $L(t)$.

A L\'{e}vy process is a real-valued stochastic phenomenon with independent and stationary increments. It illustrates the random movement of a position whose sequential displacements are erratic and unpredictable, with expulsions in pairwise separate intervals of time being independent of each other. This process can be analyzed through the continuous-time parallel of a random walk. The L\'evy process was first introduced by the French mathematician Paul L\'evy in the 1930s. The L\'evy process is particularly useful in describing systems with jumps or discontinuous movements in various fields, including probability theory, statistics, stochastic analysis, stock prices, and random turbulent flows \cite{YSD,YBD}.
Familiar examples of L\'{e}vy processes are the Pascal process, Gamma process, Poisson process, Brownian motion (also known as Wiener process), Meixner process, pure jump processes, and jump-diffusion processes \cite{ZYHD}. A L\'{e}vy process $\big\{L (t), t \geq 0\big\}$ is a stochastic process that possesses the following four main attributes:
\begin{description}
  \item[\textbf{(i)}] $L(0)=0$;
  \item[\textbf{(ii)}] Independent increments: for $t_{1}<t_{2}<\cdot\cdot\cdot<t_{m-1}<t_{m}$, the random variables $L(t_{2})-L(t_{1}),\cdot\cdot\cdot,L(t_{m})-L(t_{m-1})$ are independent;
  \item[\textbf{(iii)}]  Stationary increments: $L(t)-L(s)$ and $L(t-s)$ have the same distribution;
  \item[\textbf{(iv)}]  Stochastically continuous sampling paths (i.e., the sampling paths are continuous in probability): for all $\gamma\geq0$ and all $s\geq0$ ,
\begin{center}
$\mathbb{P}(\mid L(t)-L(s)\mid > \gamma)\rightarrow 0$, \hspace{0.4cm}\text{as}\,  $t\rightarrow s$ .
\end{center}
\end{description}

The above property \textbf{(iv)} of L\'{e}vy motion is equivalent to the following assertion: $t\mapsto L(t)$ is almost surely c\`{a}dl\`{a}g up to a modification of the process. Through equation (\ref{1}), the L\'{e}vy noise describes a method where the phase of the stimulation is perturbed.

In the recent work, generalized Kudryashov method (GKM) \cite{KBA} is utilized to attain new results in the discrete form of stochastic Schr\"odinger equation (\ref{1}). Also, the effect of L\'{e}vy noise argues on the solutions of the model \cite{PZ}. Due to our specific convention, this paper precisely examines new solutions of the stochastic nonlinear Schr\"odinger equation (\ref{1}).

In the modern efforts, our purpose is to achieve a precise understanding of the implications of the nonlinear Schr\"odinger equation (\ref{1}) through the use of stochastic processes modulated by one-dimensional L\'{e}vy noise, adopting a method such as the GKM. This approach can refine and enhance certain mathematical computations. The acquired solutions will be highly useful in elucidating various physical phenomena to physicists, including phenomena like coastal water motions, quasi-particle theory, and advancements of biomedical and fiber optics applications.

The article is decomposed into the following parts. In Section \ref{M}, we outline the methodological framework of the GKM approach, which we utilize to obtain solutions of the stochastic nonlinear Schr\"odinger equation (\ref{1}). In Section \ref{A}, we demonstrate the impact of L\'{e}vy noise on the solutions obtained for stochastic nonlinear Schr\"odinger equation, thereby  highlighting the utility of the GKM approach. In Section \ref{SA}, we present a stability analysis of the solutions. Finally, in Section \ref{CO}, we summarize our findings and provide an outlook on future research directions at the end of this paper.

\section{Methodology}\label{M}
Consider the following steps to attain the solutions:
\\\textbf{Step 1.} The general form of nonlinear PDE for a function $\psi$ of two distinct variables is
\begin{equation}\label{2}
P(\psi,\psi_{t},\psi_{x},\psi_{xx},\cdot\cdot\cdot)=0,
\end{equation}
where $P$ represents the polynomial of $\psi(x,t)$. We can get the exact solution of the travelling wave transformation of the form:
\begin{equation}\label{3}
\psi=\psi(x,t)=\psi(\xi),\hspace{0.3cm}\xi=k(x+2\alpha t).
\end{equation}
Applying the above equation (\ref{3}), the equation (\ref{2}) is converted into an ordinary differential equation (ODE):
\begin{equation}\label{4}
K(t,x,\psi,\psi',\psi'',\cdot\cdot\cdot)=0,
\end{equation}
where the subscriptions indicate the ordinary derivatives of $\psi(\xi)$ with respect to $\xi$.
\\\textbf{Step 2.} Thus, the solutions of equation (\ref{4}) will be assumed as
\begin{equation}\label{5}
\psi(\xi)=\frac{\sum_{i= 0}^{N}A_{i}\Psi^{i}(\xi)}{\sum_{j=0}^{M}B_{j}\Psi^{j}(\xi)}=\frac{T[\Psi(\xi)]}{V[\Psi(\xi)]},
\end{equation}
where $A_{i} (i=0,1,\cdot\cdot\cdot,N)$ and $B_{j} (j=0,1,\cdot\cdot\cdot,M)$ are parameters to be analyzed such that $A_{N}\neq0$ and $B_{M}\neq0$, while $\Psi(\xi)$ has the set up
\begin{equation}\label{6}
\Psi(\xi) = \frac{1}{1 + Ae^{\xi}},
\end{equation}
where $A$ is arbitrary constant. The solution to the Ricatti equation is
\begin{equation*}
\Psi'(\xi)=\Psi^{2}(\xi)-\Psi(\xi).
\end{equation*}
Utilizing equation (\ref{6}), the following derivatives are attained:
\begin{eqnarray}\label{7}
\psi'(\xi)&=&\frac{T'\Psi'V-TV'\Psi'}{V^{2}}=\Psi'\big[\frac{T'V-TV'}{V^{2}}\big]=(\Psi^{2}-\Psi)\big[\frac{T'V-TV'}{V^{2}}\big],\\ \label{8}
\psi''(\xi)&=&\frac{\Psi^{2}-\Psi}{V^{2}}\big[(2\Psi-1)(T'V-TV')\big]+\frac{\Psi^{2}-\Psi}{V^{3}}\big[V(T''V-TV'')-2V'T'V+2T(V')^{2}\big].
\end{eqnarray}
\\\textbf{Step 3.} The desired solution of the nonlinear ODE obtained by equation (\ref{4}) employing the GKM approach is:
\begin{equation}\label{9}
\psi(\xi) = \frac{A_{0}+A_{1}\Psi+A_{2}\Psi^{2}+\cdot\cdot\cdot+A_{N}\Psi^{N}}{B_{0}+B_{1}\Psi+B_{2}\Psi^{2}+\cdot\cdot\cdot+B_{M}\Psi^{M}}.
\end{equation}
Apply the homogeneous balance principle to evaluate the values of $M$ and $N$ in equation (\ref{5}) by balancing the highest-order derivative and nonlinear term in equation (\ref{4}).
\\\textbf{Step 4.} Substitute equation (\ref{5}) into equation (\ref{4}) and thus obtain a polynomial $K(\Psi(\xi)$) of $\Psi$. Furthermore, by equating each coefficient of $P(\Psi(\xi)$) to zero, we obtain an algebraic system. We then solve this system to determine the values of the unknown coefficients $A_{0},A_{1},A_{2},\cdot\cdot\cdot,A_{N},B_{0},B_{1},B_{2},\cdot\cdot\cdot,B_{M}$. Finally, we analyze the solutions of nonlinear PDE (\ref{2}) obtained through the previous steps.

\section{Applications}\label{A}
In this section, we aim to attain the exact solution of  the stochastic nonlinear Schr\"odinger equation (\ref{1}). Consider the traveling wave transformation:
\begin{equation}\label{10}
\psi(x,t)=e^{i\theta}u(\xi),\quad \xi=k(x+2\alpha t),\quad \theta=\alpha x+\upsilon t+\sigma L(t),
\end{equation}
where $k$ is a positive constant, $\sigma$ is the intensity of L\'{e}vy noise, and $\alpha$ is the wave speed of $u(\xi)$. Through equation (\ref{10}), we obtain the following expressions:
\begin{align}\label{11}
\begin{split}
 \psi_{xx}&=e^{i\theta}(k^2u''+2i\alpha ku'-\alpha^2u),
\\
 \psi_{t}&=e^{i\theta}(2\alpha ku'+i\upsilon u+i\sigma uL_{t}).
\end{split}
\end{align}
By substituting these expressions into equation (\ref{1}) and rearranging, we arrive at
the following ODE representation:
\begin{equation}\label{12}
-k^2u''+2u^3+Hu=0.
\end{equation}
Thus,
\begin{equation*}
H=\alpha^2-2\rho^2-\upsilon.
\end{equation*}
By applying the homogenous balance principle to the ODE equation (\ref{12}), we attain a positive integer $N=M+1$ by balancing the terms $u''$ and $u^3$. We choose integers $M=1$ and $N=2$ to obtain the following specific solution of the form:
\begin{equation}\label{13}
u(\xi)=\frac{A _{0}+A _{1}\Psi+A _{2}\Psi^{2}}{B_{0}+B_{1}\Psi},
\end{equation}
and therefore, we can express the derivative as:
\begin{equation}\label{14}
u'(\xi)=(\Psi^{2}-\Psi)\Big[\frac{(A_{1}+2A_{2}\Psi)(B_{0}+B_{1}\Psi)-B_{1}(A_{0}+A_{1}\Psi+A_{2}\Psi^{2}}{(B_{0}+B_{1}\Psi)^{2}}\Big],
\end{equation}
\begin{equation}\label{15}
\begin{aligned}
u''(\xi) &= \frac{\Psi^{2}-\Psi}{(B_{0}+B_{1}\Psi)^{2}}(2\Psi-1)\big[(A_{1}+2A_{2}\Psi)(B_{0}+B_{1}\Psi)-B_{1}(A_{0}+A_{1}\Psi+A_{2}\Psi^{2}\big]\\
&\,\,\,\,\,+\frac{(\Psi^{2}-\Psi)^{2}}{(B_{0}+B_{1}\Psi)^{3}}\big[2A_{2}(B_{0}+B_{1}\Psi)^{2}-2B_{1}(A_{1}+2A_{2}\Psi)(B_{0}+B_{1}\Psi)+2B_{1}^{2}(A_{0}+A_{1}\Psi+A_{2}\Psi^{2})\big].
\end{aligned}
\end{equation}
We obtain the polynomial in $\Psi(\xi)$ by considering equation (\ref{13}) and shifting equation (\ref{15}) along with its required derivatives. Then we collect the variable coefficients corresponding to the same power terms of
$\Psi(\xi)$ and set them equal to zero, thereby obtaining an algebraic system of equations:
\begin{equation}\label{16}
\begin{aligned}
\Psi^{0}(\xi)&: 2A_{0}^{3} + H A_{0}B_{0}^{2} = 0,\\
\Psi^{1}(\xi)&: 6A_{0}^{2}A_{1} + H A_{1}B_{0}^{2} - k^{2} A_{1} B_{0}^{2} + 2H A_{0} B_{0} B_{1} + k^{2} A_{0} B_{0} B_{1} = 0,\\
\Psi^{2}(\xi)&: 6A_{0} A_{1}^{2} + 6A_{0}^{2} A_{2} + 3k^{2} A_{1} B_{0}^{2} - 4k^{2} A_{2} B_{0}^{2} - 3k^{2} A_{0} B_{0} B_{1}\\ &+ 2H A_{1} B_{0} B_{1} + k^{2} A_{1} B_{0} B_{1} + H A_{0} B_{1}^{2} - k^{2} A_{0} B_{1}^{2} = 0,\\
\Psi^{3}(\xi)&: 2A_{1}^{3} + 12A_{0} A_{1} A_{2} -2k^{2} A_{1} B_{0}^{2} + 10k^{2} A_{2} B_{0}^{2} + 2k^{2} A_{0} B_{0} B_{1}\\ &-k^{2} A_{1} B_{0} B_{1}+ 2H A_{2} B_{0} B_{1} -3k^{2} A_{2} B_{0} B_{1} +k^{2} A_{0} B_{1}^{2} + H A_{1} B_{1}^{2} = 0,\\
\Psi^{4}(\xi)&: 6A_{1}^{2} A_{2} +6A_{0} A_{2}^{2} -6k^{2} A_{2} B_{0}^{2} + 9k^{2} A_{2} B_{0} B_{1} + H A_{2} B_{1}^{2} -k^{2} A_{2} B_{1}^{2} = 0,\\
\Psi^{5}(\xi)&: 6A_{1}  A_{2}^{2} -6k^{2} A_{2} B_{0} B_{1} + 3k^{2} A_{2} B_{1}^{2}= 0,\\
\Psi^{6}(\xi)&: 2A_{2}^{3} - 2k^{2} A_{2} B_{1}^{2}= 0.
\end{aligned}
\end{equation}
Subsequently, we analyze the above algebraic system of equations to determine the values of the unknown arbitrary constants $A_{0}$, $A_{1}$, $A_{2}$, $B_{0}$, $B_{1}$ and $k$. After establishing the solution sets of the system, we can obtain the solutions for the nonlinear Schr\"odinger equation (\ref{1}) in the following specific cases.
\\
\\
$\textbf{Case 1}$:
\begin{eqnarray*}
k=i\sqrt{2H}, \hspace{0.2cm}A_{0}=\frac{i\sqrt{H} B_{0}}{\sqrt{2}}, \hspace{0.2cm}A_{1}=-i\sqrt{2H} B_{0},\hspace{0.2cm}A_{2}= 0, \hspace{0.2cm}B_{1} =-2B_{0}.
\end{eqnarray*}
The solution of equation (\ref{1}) is obtained by setting the case 1 in equation (\ref{13}) as:
\begin{eqnarray}\label{17}
\Psi_{1}(x,t) = \left(\frac{\frac{i\sqrt{H} B_{0}}{\sqrt{2}}-\frac{i\sqrt{2H} B_{0}}{1+A e^{i\sqrt{2H} (x+2\alpha t)}}}                {B_{0}-\frac{2B_{0}}{1+A e^{i\sqrt{2H} (x+2\alpha t)}}}\right)e^{i\big(\alpha x+\upsilon t+\sigma L(t)\big)}.
\end{eqnarray}
\\
$\textbf{Case 2}$:
\begin{eqnarray*}
k=i\sqrt{2H}, \hspace{0.2cm}A_{0}=\frac{i\sqrt{H} B_{0}}{\sqrt{2}}, \hspace{0.2cm}A_{1}=0, \hspace{0.2cm} A_{2}=0, \hspace{0.2cm}B_{1}=-2B_{0}.
\end{eqnarray*}
The solution of equation (\ref{1}) is derived by setting the case 2 in equation (\ref{13}) as:
\begin{eqnarray}\label{18}
\Psi_{2}(x,t)=\left(\frac{i\sqrt{H} B_{0}}{\sqrt{2}(B_{0}-\frac{2B_{0}}{1+A e^{i\sqrt{2H}(x+2\alpha t)}})}\right)e^{i\big(\alpha x+\upsilon t+\sigma L(t)\big)}.
\end{eqnarray}
\\
$\textbf{Case 3}$:
\begin{eqnarray*}
k=i\sqrt{2H}, \hspace{0.2cm}A_{0}=-\frac{i\sqrt{H}B_{0}}{\sqrt{2}}, \hspace{0.2cm}A_{1}=i\sqrt{2H}B_{0}, \hspace{0.2cm}A_{2}=0, \hspace{0.2cm}B_{1} =-2B_{0}.
\end{eqnarray*}
The solution of equation (\ref{1}) is gotten by setting the case 3 in equation (\ref{13}) as:
\begin{eqnarray}\label{19}
\Psi_{3}(x,t)=\left(\frac{-\frac{i\sqrt{H} B_{0}}{\sqrt{2}}+\frac{i\sqrt{2H} B_{0}}{1+A e^{i\sqrt{2H} (x+2\alpha t)}}}                {B_{0}-\frac{2B_{0}}{1+A e^{i\sqrt{2H} (x+2\alpha t)}}}\right)e^{i\big(\alpha x+\upsilon t+\sigma L(t)\big)}.
\end{eqnarray}
\\
$\textbf{Case 4}$:
\begin{eqnarray*}
k=-\frac{i\sqrt{H}}{\sqrt{2}}, \hspace{0.2cm}A_{0}=-\frac{i\sqrt{H}B_{0}}{\sqrt{2}}, \hspace{0.2cm}A_{1}=i\sqrt{2H} B_{0}, \hspace{0.2cm}A_{2}= -i\sqrt{2H}B_{0}, \hspace{0.2cm}B_{1}=-2B_{0}.
\end{eqnarray*}
The solution of equation (\ref{1}) is gained by setting the case 4 in equation (\ref{13}) as:
\begin{eqnarray}\label{20}
\Psi_{4}(x,t)=\left(\frac{-\frac{i\sqrt{H}B_{0}}{\sqrt{2}}-\frac{i\sqrt{2H}B_{0}}{(1+A e^{-i\sqrt{H}(x+2\alpha t)})^{2}}+\frac{i\sqrt{2H}B_{0}}{1+A e^{-i\sqrt{H}(x+2\alpha t)}}}{B_{0}-\frac{2B_{0}}{1+Ae^{-i\sqrt{H}(x+2\alpha t)}}}\right)e^{i\big(\alpha x+\upsilon t+\sigma L(t)\big)}.
\end{eqnarray}
\\
$\textbf{Case 5}$:
\begin{eqnarray*}
k=\frac{i\sqrt{H}}{\sqrt{2}}, \hspace{0.2cm}A_{0}=\frac{i\sqrt{H} B_{0}}{\sqrt{2}}, \hspace{0.2cm}A_{1}=-i\sqrt{2H} B_{0}, \hspace{0.2cm}A_{2}=i \sqrt{2H}B_{0}, \hspace{0.2cm}B_{1}=-2B_{0}.
\end{eqnarray*}
The solution of equation (\ref{1}) is acquired by setting the case 5 in equation (\ref{13}) as:
\begin{eqnarray}\label{21}
\Psi_{5}(x,t)=\left(\frac{\frac{i\sqrt{H}B_{0}}{\sqrt{2}}+\frac{i\sqrt{2H}B_{0}}{(1+Ae^{-i\sqrt{H}(x+2\alpha t)})^{2}}-\frac{i\sqrt{2H}B_{0}}{1+A e^{-i\sqrt{H}(x+2\alpha t)}}}{B_{0}-\frac{2B_{0}}{1+A e^{-i\sqrt{H}(x+2\alpha t)}}}\right)e^{i\big(\alpha x+\upsilon t+\sigma L(t)\big)}.
\end{eqnarray}
\\
$\textbf{Case 6}$:
\begin{eqnarray*}
k=-i\sqrt{2H}, \hspace{0.2cm}A_{0}=-\frac{i\sqrt{H} B_{0}}{\sqrt{2}}, \hspace{0.2cm}A_{1}=0, \hspace{0.2cm}A_{2}=0, \hspace{0.2cm}B_{1}=-2B_{0}.
\end{eqnarray*}
The solution of equation (\ref{1}) is achieved by setting the case 6 in equation (\ref{13}) as:
\begin{eqnarray}\label{22}
\Psi_{6}(x,t)=\left(-\frac{i\sqrt{H}B_{0}}{\sqrt{2}(B_{0}-\frac{2B_{0}}{1+Ae^{-i\sqrt{2H}(x+2\alpha t)}})}\right)e^{i\big(\alpha x+\upsilon t+\sigma L(t)\big)}.
\end{eqnarray}
\\
$\textbf{Case 7}$:
\begin{eqnarray*}
k=-i\sqrt{2H}, \hspace{0.2cm}A_{0}=\frac{i\sqrt{H} B_{0}}{\sqrt{2}}, \hspace{0.2cm}A_{1}=-i\sqrt{2H}B_{0}, \hspace{0.2cm}A_{2}=0, \hspace{0.2cm}B_{1} =-2B_{0}.
\end{eqnarray*}
The solution of equation (\ref{1}) is obtained by setting the case 7 in equation (\ref{13}) as:
\begin{eqnarray}\label{23}
\Psi_{7}(x,t)=\left(\frac{\frac{i\sqrt{H} B_{0}}{\sqrt{2}}-\frac{i\sqrt{2H} B_{0}}{1+Ae^{-i\sqrt{2H}(x+2\alpha t)}}}                {B_{0}-\frac{2B_{0}}{1+A e^{-i\sqrt{2H}(x+2\alpha t)}}}\right)e^{i\big(\alpha x+\upsilon t+\sigma L(t)\big)}.
\end{eqnarray}
\\
$\textbf{Case 8}$:
\begin{eqnarray*}
A_{0}=\frac{i\sqrt{2H}(2B_{0}^{2}+B_{0}B_{1})}{2(2B_{0}+B_{1})}, \hspace{0.2cm}A_{1}=\frac{i\sqrt{H}B_{0}}{\sqrt{2}}, \hspace{0.2cm}A_{2}=0.
\end{eqnarray*}
The solution of equation (\ref{1}) is procured by setting the case 8 in equation (\ref{13}) as:
\begin{eqnarray}\label{24}
\Psi_{8}(x,t)=\left(\frac{\frac{i\sqrt{H}B_{1}}{\sqrt{2}(1+Ae^{k(x+2\alpha t)})}+\frac{i\sqrt{H}(2B_{0}^{2}+B_{0}B_{1})}{\sqrt{2}(2B_{0}+B_{1})}}{B_{0}+\frac{B_{1}}{1+Ae^{k(x+2\alpha t)}}}\right)e^{i\big(\alpha x+\upsilon t+\sigma L(t)\big)}.
\end{eqnarray}
\\
$\textbf{Case 9}$:
\begin{eqnarray*}
A_{0}=-\frac{i\sqrt{2H}(2B_{0}^{2}+B_{0}B_{1})}{2(2B_{0}+B_{1})}, \hspace{0.2cm}A_{1}=-\frac{i\sqrt{H}B_{0}}{\sqrt{2}}, \hspace{0.2cm}A_{2}=0.
\end{eqnarray*}
The solution of equation (\ref{1}) is gotten by setting the case 9 in equation (\ref{13}) as:
\begin{eqnarray}\label{25}
\Psi_{9}(x,t)=-\left(\frac{\frac{i\sqrt{H}B_{1}}{\sqrt{2}(1+Ae^{k(x+2\alpha t)})}+\frac{i\sqrt{H}(2B_{0}^{2}+B_{0}B_{1})}{\sqrt{2}(2B_{0}+B_{1})}}     {B_{0}+\frac{B_{1}}{1+Ae^{k(x+2\alpha t)}}}\right)e^{i\big(\alpha x+\upsilon t+\sigma L(t)\big)}.
\end{eqnarray}
\\
$\textbf{Case 10}$:
\begin{eqnarray*}
k=-i\sqrt{2H}, \hspace{0.2cm}A_{0}=\frac{i\sqrt{H} B_{0}}{\sqrt{2}}, \hspace{0.2cm}A_{1}=i\sqrt{2H}(-B_{0}+B_{1}), \hspace{0.2cm}A_{2}=-i\sqrt{2H}B_{1}.
\end{eqnarray*}
The solution of equation (\ref{1}) is gained by setting the case 10 in equation (\ref{13}) as:
\begin{eqnarray}\label{26}
\Psi_{10}(x,t)=-\left(\frac{\frac{i\sqrt{H}B_{1}}{\sqrt{2}(1+Ae^{k(x+2\alpha t)})}+\frac{i\sqrt{H}(2B_{0}^{2}+B_{0}B_{1})}{\sqrt{2}(2B_{0}+B_{1})}}     {B_{0}+\frac{B_{1}}{1+Ae^{k(x+2\alpha t)}}}\right)e^{i\big(\alpha x+\upsilon t+\sigma L(t)\big)}.
\end{eqnarray}
\\
$\textbf{Case 11}$:
\begin{eqnarray*}
k=-i\sqrt{2H}, \hspace{0.2cm}A_{0}=\frac{i\sqrt{2H}(2B_{0}^{2}+B_{0}B_{1})}{2(2B_{0}+B_{1})}, \hspace{0.2cm}A_{1}=-\frac{ i\sqrt{H}(2B_{0}+B_{1})}{\sqrt{2}}, \hspace{0.2cm}A_{2}=0.
\end{eqnarray*}
The solution of equation (\ref{1}) is acquired by setting the case 11 in equation (\ref{13}) as:
\begin{eqnarray}\label{27}
\Psi_{11}(x,t)=\left(\frac{-\frac{i\sqrt{2H}B_{1}}{\sqrt{2}(1+Ae^{-i\sqrt{2H}(x+2\alpha t)})}+\frac{i\sqrt{H}(2B_{0}^{2}+B_{0}B_{1})}{\sqrt{2}(2B_{0}+B_{1})}}{B_{0}+\frac{B_{1}}{1+Ae^{-i\sqrt{2H}(x+2\alpha t)}}}\right) \hspace{0.2cm}e^{i\big(\alpha x+\upsilon t+\sigma L(t)\big)}.
\end{eqnarray}
\\
$\textbf{Case 12}$:
\begin{eqnarray*}
k=i\sqrt{2H}, \hspace{0.2cm}A_{0}=\frac{i\sqrt{2H}(2B_{0}^{2}+B_{0}B_{1})}{2(2B_{0}+B_{1})}, \hspace{0.2cm}A_{1}=\frac{ i\sqrt{H}(2B_{0}+B_{1})}{\sqrt{2}}, \hspace{0.2cm}A_{2}=0.
\end{eqnarray*}
The solution of equation (\ref{1}) is achieved by setting the case 12 in equation (\ref{13}) as:
\begin{eqnarray}\label{28}
\Psi_{12}(x,t)=\left(\frac{\frac{i\sqrt{2H}B_{1}}{\sqrt{2}(1+Ae^{-i\sqrt{2H}(x+2\alpha t)})}+\frac{i\sqrt{H}(2B_{0}^{2}+B_{0}B_{1})}{\sqrt{2}(2B_{0}+B_{1})}}{B_{0}+\frac{B_{1}}{1+Ae^{-i\sqrt{2H}(x+2\alpha t)}}}\right)e^{i\big(\alpha x+\upsilon t+\sigma L(t)\big)}.
\end{eqnarray}
Therefore, the solutions of the nonlinear Schr\"odinger equation (\ref{1}) with L\'{e}vy noise were found precisely by utilizing the generalized Kudryashov method.

\section{Stability Analysis}\label{SA}
Because of the analysis of the NLS (\ref{1}), which is based on the stability arising from a specific form of the nonlinearity \cite{36}, the momentum for equation (\ref{1}) can be described as
\begin{equation}\label{29}
Q=\frac{1}{2}\int_{-\infty}^{\infty}\nu^{2}(\xi)d\xi,
\end{equation}
where $\nu$ is the electric potential and $Q$ is the momentum. The stability or soliton state of the solitary wave for this model is given by
\begin{equation}\label{30}
\frac{\partial Q}{\partial \lambda}>0,
\end{equation}
where $\lambda$ is the frequency. By substituting the solitary wave equation (\ref{18}) into equation (\ref{29}), we acquire the following structure
\begin{equation}\label{31}
Q=\int_{-10}^{10}\kappa^{2}\left(\frac{i\sqrt{H} B_{0}}{\sqrt{2}(B_{0}-\frac{2B_{0}}{1+Ae^{i\sqrt{2H}(x+2\alpha t)}})}\right)d\xi.
\end{equation}
By evaluating the integral in equation (\ref{31}), we arrive at the following expression:
\begin{equation*}
Q=-\frac{\kappa^{2}}{2}H\left(20-\frac{i2\sqrt{2}(\frac{1}{-1+Ae^{i2\sqrt{2H}(-5+2\alpha)}}-\frac{1}{-1+Ae^{i2\sqrt{2H}(5+2\alpha)}})} {\sqrt{H}}\right).
\end{equation*}
Now, we apply the stability condition expressed in equation (\ref{30}) to the above expression and thus get
\begin{align*}
&\frac{\kappa^{2}}{(-1+Ae^{i4(-5+2\alpha)})^{2}(-1+Ae^{i4(5+2\alpha)})^{2}}\Big[-10+\frac{1}{2}i\Big(\frac{1}{-1+Ae^{i4(-5+2\alpha)}}+\frac{1}{1+Ae^{i4(5+2\alpha)}}\Big)\\&-2A\Big(-20Ae^{i16\alpha}+A^{2}e^{i4(-5+6\alpha)}\big(5+e^{i40}(5-2\alpha)+2\alpha\big)+e^{i4(-5+2\alpha)}\big(5-2\alpha+e^{i40}(5+2\alpha)\big)\Big) \Big]>0.
\end{align*}
Therefore, it indicates that equation (\ref{1}) is a stable nonlinear equation.

Similarly, by substituting the wave equation (\ref{22}) into equation (\ref{30}), we obtain the following form:
\begin{eqnarray}\label{32}
Q=\int_{-10}^{10}\kappa^{2}\left(-\frac{i\sqrt{H} B_{0}}{\sqrt{2}(B_{0}-\frac{2B_{0}}{1+Ae^{i\sqrt{2H}(x+2\alpha t)}})}\right)d\xi.
\end{eqnarray}
The following expression is derived by evaluating the above integral:
\begin{eqnarray*}
Q=\kappa^{2}\left(-\frac{i\sqrt{2H}}{-1+Ae^{i2\sqrt{2H}(5-2\alpha)}}-\sqrt{2}i(1+\frac{A}{-A+e^{i2\sqrt{2H}(5+2\alpha)}})\sqrt{H}-10H\right).
\end{eqnarray*}
Now, the following analysis is applied. By utilizing the stability state expressed in equation (\ref{30}) to above equation, we have
\begin{eqnarray*}
\kappa^{2}\Big[-10-\frac{1}{2}i\big(1+\frac{1}{-1+Ae^{i4(5-2\alpha)}}\big)+A\big(\frac{i}{2A-2e^{i4(5+2\alpha)}}+\frac{2e^{i4(5+2\alpha)}(-5+2\alpha)}{(-Ae^{i20}+e^{i8\alpha})^{2}}-\frac{2e^{i4(5+2\alpha)}(5+2\alpha)}{(A-e^{i4(5+2\alpha)})^{2}}\big)\Big]>0.
\end{eqnarray*}
Hence,  equation (\ref{1}) is also a stable nonlinear equation.

\section{Conclusion and Outlook}\label{CO}

In this paper, we have introduced the methodological backbone of the generalized Kudryashov method, which serves as the primary tool for deriving solutions to the stochastic nonlinear Schr\"{o}dinger equation \eqref{1}. Through the application of the generalized Kudryashov method, we have identified a more diverse set of stochastic solutions under specific conditions. Additionally, we have investigated how L\'evy noise impacts the solutions of the nonlinear Schr\"{o}dinger equation and interpreted the dynamic characteristics of the results. Furthermore, we have executed a thorough stability analysis of the derived solutions.

In the future, we will continue the research in the area of stochastic nonlinear Schr\"{o}dinger equations, focusing on various methodologies, approaches, and techniques for solving specific problems in the context of L\'evy noise's influence.  We expect to see more interdisciplinary collaborations between different research communities to address complex problems related to L\'evy noise and the generalized Kudryashov method. Furthermore, we anticipate further development and application of solitary wave configurations in physical and engineering areas such as nonlinear optics, plasma physics, Bose-Einstein condensates, and water waves. Finally, we believe that the findings and contributions of this paper will provide valuable insights and guidance for future research in the field of novel stochastic methods, specifically for nonlinear Schr\"{o}dinger equations, that are more efficient, accurate, and reliable.


\end{document}